\documentclass[12pt,reqno]{article}
\usepackage{dcolumn}
\usepackage{amsthm,epsfig}
\usepackage{array,dcolumn,graphicx,amsmath}
\usepackage{booktabs,lscape,multirow,psfrag,rotating}
\usepackage{amssymb,amsmath,amsthm}
\usepackage{longtable}
\usepackage{booktabs,lscape,multirow,rotating}
\usepackage[normal]{caption2}
\usepackage[mathscr]{eucal}
\usepackage{eucal}
\usepackage{amssymb,pstricks}
\newcommand{\bSigma}{\boldmath$\Sigma$}
\newcommand{\mbSigma}{\mbox{\bSigma}}

\newcommand{\bGamma}{\boldmath$\Gamma$}
\newcommand{\mbGamma}{\mbox{\bGamma}}

\newcommand{\btau}{\boldmath$\tau$}
\newcommand{\mbtau}{\mbox{\btau}}

\newcommand{\bPhi}{\boldmath$\Phi$}
\newcommand{\mbPhi}{\mbox{\bPhi}}
\newcommand{\bPsi}{\boldmath$\Psi$}
\newcommand{\mbPsi}{\mbox{\bPsi}}

\newcolumntype{d}[1]{D{.}{.}{#1}}
\begin{document}
\def\s#1{\oalign{$#1$\crcr\hidewidth \normal$\sim$ \hidewidth}}
\thispagestyle{empty}
\newtheorem{theorem}{\indent THEOREM}
\newtheorem{prop}{\indent Proposition}
\newtheorem{lemma}{\indent Lemma}
\renewcommand{\proofname}{\hspace*{\parindent}{Proof.}}

\begin{center}
{\bf On the orthogonally equivariant estimators of a covariance
matrix}
\end{center}

\begin{center}
{Ming-Tien Tsai and Chia-Hsuan Tsai}
\end{center}

\begin{center}
{\it Institute of Statistical Science, Academia Sinica, Taipei, Taiwan 11529, R.O.C.}
\end{center}

\vspace{0.3cm}
\begin{center}
\parbox{14cm}{\sloppy \, \,
{\bf Abstract}. In this note, when the dimension $p$ is large we look into the insight of the Mar$\check{c}$enko-Pastur
equation to get an explicit equality relationship, and use the obtained equality to establish a new orthogonally
equivariant estimator of the population covariance matrix. Under some regularity conditions, the proposed novel estimators
of the population eigenvalues are shown to be consistent with the eigenvalues of the population covariance matrix. It is
also demonstrated that the proposed estimator is the best orthogonally equivariant estimator of the population covariance
matrix under the normalized Stein loss function.} \\

\vspace{0.3cm}
\small
\parbox{14cm}{\sloppy \, \,
{\bf Keywords:} Consistent estimators of the population eigenvalues, orthogonally equivariant estimators,
high-dimensional covariance matrix.}
\small
\parbox{14cm} {\small
}
\end{center}

\vspace{0.5cm}
\def \theequation{1.\arabic{equation}}
\setcounter{equation}{0}
\noindent {\bf 1. Introduction}
\vspace{0.3cm}

\indent The problem in high-dimensional covariance estimation has been one of the most interesting topics in statistics
(Pourahmadi, 2013; Zagidullina, 2021). Stein (1975, 1986) considered the orthogonally equivariant nonlinear shrinkage
estimator for the population covariance matrix. Stein's estimator has been considered a gold standard, and from which
a large strand of literature on the orthogonally equivariant estimation of covariance matrix was generated (Ledoit and Wolf,
2012; Rajaratnam and Vincenzi, 2016, and the references therein). Under the large dimensional asymptotics setup, namely,
both sample size $n$ and the dimension $p$ are sufficiently large with the concentration
$c=\lim_{n \to \infty}{p}/{n}, c \in (0, \infty)$, random matrix theory has been incorporated to find out the reasonable
estimators of the population covariance matrix. Ledoit and Wolf (2012) used the generalized Mar$\check{c}$enko-Pastur
equation of Ledoit and P$\acute{e}$ch$\acute{e}$ (2011) to propose their orthogonally equivariant estimators for the
population covariance matrix and its precision matrix.

\indent Ledoit and Wolf (2018) commented that Stein's estimator has its theoretical limitations, and further they asserted
the asymptotic optimality of their estimators under the normalized Stein loss function as well as some other modified loss
functions. The difference between the two estimators is to use different ways to estimate the Stieltjes transform of limiting
empirical spectral distribution function. Stein (1986) used the naive empirical distribution of sample eigenvalues to
estimate the Stieltjes transform of the distribution function, while Ledoit and Wolf (2012) used a more sophisticated
smoothed version called the QuEST function instead. Both those two estimators of Stein (1975, 1986) and Ledoit and Wolf
(2012) are all consistent for the Stieltjes transform of limiting empirical spectral distribution function, hence
theoretically they should enjoy the same asymptotical property

\vspace{0.3cm}
\noindent \underline{~~~~~~~~~~~~~~~~~~~~~~~~~~~~~~~~~~~~~~~~~~~~}\\
E-mail: mttsai@stat.sinica.edu.tw
\pagebreak

\noindent in the sense of large dimensional asymptotics although Ledoit and Wolf's estimator may have some better
numerical performances in practice. It is easy to see that both Stein's estimator and Ledoit and Wolf's estimator are
asymptotically equivalent under the large dimensional asymptotics setup, and all are reduced to the sample covariance
matrix when the dimension $p$ is fixed and the sample size $n$ is large.

\indent Although the Mar$\check{c}$enko-Pastur equation (Silverstein, 1995) provides the link of limiting empirical
spectral distribution, $F$, to the limiting behavior of the population spectral distribution, $H$, one can expect to
retrieve the information of $H$ from $F$. However, the difficulty lies in the fact that the relationship between
$H$ and $F$ are entangled (see (4.3)). This phenomenon has been open for more than 50 years in the literature since 1967.

\indent Ledoit and Wolf (2012) generalized the Mar$\check{c}$enko-Pastur equation to design for the population covariance
matrix $\mbSigma$ and the precision matrix $\mbSigma^{-1}$, respectively. However, the product of two obtained estimators
by their procedure does not converge to the identity matrix ${\bf I}$ in probability, hence their results seem to be
against statistical senses. Thus, it is reasonable to suspect that their estimator of $\mbSigma$ might not be optimal among
the class of orthogonally equivalent estimators. Therefore, there remains room for further investigation of the estimation
problem of $\mbSigma$.

\indent As such, the other direction is to see whether another type of explicit expression of the relationship of $F$ and
$H$ can be further dug out. When the dimension $p$ is large, we look into the insight of the Mar$\check{c}$enko-Pastur
equation, by mimicking the proof of the Sokhotski-Plemelj formula, and luckily get another kind of explicit equality
relationship between the quantiles of $F$ and $H$.

\indent Some notations of the orthogonally equivariant estimator are introduced in Section 2. In the same section when the
dimension $p$ is fixed we claim that the sample covariance matrix is the best orthogonally equivariant estimator under the
Stein loss function. In Section 3, when the dimension $p$ is fixed, we show that Stein's estimator can be inadmissible
under the Stein loss function. In Section 4, when the dimension $p$ is large under the large dimensional asymptotics setup,
namely both sample size $n$ and the dimension $p$ are sufficiently large with the concentration
$\lim_{n \to \infty} p/n=c \in [0, 1)$, we re-examine the asymptotic optimal property of estimators claimed by both
Stein (1975, 1986) and Ledoit and Wolf (2018). Moreover, the obtained new kind of explicit equalities between the quantiles
of $F$ and $H$ are used to establish a new type of orthogonally equivariant estimator for the population covariance matrix.
We show that our proposed estimators of the eigenvalues are consistent with the eigenvalues of the population covariance
matrix. We also show that under the normalized Stein loss function the proposed estimator is the best orthogonally
equivariant estimator of the population covariance matrix. In contrast, Stein's estimator and the sample covariance
matrix are inadmissible.

\vspace{0.3cm}
\noindent {\bf 2. Preliminary notations}
\def \theequation{2.\arabic{equation}}
\setcounter{equation}{0}
\vspace{0.3cm}

\indent  Let ${\bf X}_{1}, \ldots, {\bf X}_{n}$ be independent $p$-dimensional random vectors with a common
multivariate normal distribution $N_{p}(\bf 0, {\mbSigma})$. A basic problem considered in the literature is the
estimation of the $p \times p$ covariance matrix ${\mbSigma}$, which is unknown and assumed to be non-singular.
It is also assumed that $n > p$, as such the sufficient statistic
\begin{align}
{\bf A}=\sum_{i=1}^{n} {\bf X}_{i}{\bf X}^{\top}_{i}
\end{align}
is positive definite with probability one. In the literature, the estimators $\phi({\bf A})$ of ${\mbSigma}$ are the
functions of ${\bf A}$. The sample space ${\mathcal S}$, the parameter space ${\varTheta}$, and the action space
${\mathcal A}$ are taken to be the set $\mathcal P_{p}$ of $p \times p$ symmetric positive definite matrices. Note
that ${\bf A}$ has a Wishart distribution $W({\mbSigma}, n)$, and the maximum likelihood estimator (MLE) of ${\mbSigma}$
is expressed as below
\begin{align}
\Hat{\mbSigma}_{ML}={\bf S}, ~\mbox {where} ~{\bf S}=n^{-1}{\bf A}
\end{align}
which is unbiased (Anderson, 2003). The general linear group $Gl(p)$ acts on the space $\mathcal P_{p}$.

\indent One of the most interesting loss functions was introduced by Stein (1956)
\begin{align}
L(\phi({\bf S}), {\mbSigma})=\mbox{tr}{\mbSigma}^{-1}\phi({\bf S})
      -\mbox{log}\mbox{det}{\mbSigma}^{-1}\phi({\bf S})-p,
\end{align}
where $\mbox{tr}$ and $\mbox{det}$ denote the trace and the determinant of a matrix, respectively. Because
$Gl(p)$ acts transitively on the space $\mathcal P_{p}$, the best $Gl(p)$-equivariant estimator exists.

\indent We consider invariant loss function $L$, i.e., $L$ satisfies the condition that
$L(g\phi({\bf S})g^{\top}, g{\mbSigma}g^{\top})\linebreak=L(\phi({\bf S}), {\mbSigma})$ for all $g \in Gl(p)$.
An estimator $\Hat{\mbSigma}$ is called $\mathcal O(p)$-equivariant if $\Hat{\mbSigma}({\bf G}{\bf S}{\bf G}^{\top})
={\bf G} \Hat{\mbSigma}({\bf S}){\bf G}^{'}, \forall {\bf G} \in \mathcal O(p), \forall {\bf S} \in \mathcal P_{p}$,
where $\mathcal O(p)$ is the set of orthogonal groups. Suppose that ${\bf G}$ acts on $\mathcal P_{p}$, whereby the
orbit through $x \in \mathcal P_{p}$ is the set ${\bf G}x=\{gx|g \in {\bf G}\} \subset \mathcal P_{p}$. The action is
called transitive if ${\varTheta}$ is one orbit, i.e., $\forall x, y \in {\varTheta}$ there is some $g \in {\bf G}$
with $gx=y$. It may then be easy to note the fact that if $L$ is $\mathcal O(p)$-invariant, $\Hat{\mbSigma}$ is
$\mathcal O(p)$-equivariant, and ${\bf G}$ acts transitively on $\mathcal P_{p}$, then the risk function is constant
on the $\mathcal O(p)$-orbits of $\mathcal P_{p}$:
\begin{align}
R(\Hat {\mbSigma}, {\mbSigma})=R(\Hat {\mbSigma}, {\mbGamma}),
\end{align}
where ${\mbGamma}$ is a diagonal matrix of eigenvalues of ${\mbSigma}, \forall {\mbSigma} \in \mathcal P_{p}$ and
$R(\Hat {\mbSigma}, {\mbSigma})={\mathcal E}L(\Hat {\mbSigma}, {\mbSigma})$ with ${\mathcal E}({\bf X})$ denoting
the expectation of $X$.

\indent The class of orthogonally equivariant (i.e., $\mathcal O(p)$-equivariant) estimators of covariance matrix is
constituted of all the estimators having the same eigenvectors as the sample covariance matrix. Consider the
spectral decomposition of population matrix, namely ${\mbSigma}={\bf V}{\mbGamma}{\bf V}^{\top}$, where
${\mbGamma}$ is a diagonal matrix with eigenvalues $\gamma_{i, p}, i=1, \ldots, p$, and
${\bf V}=( {\bf v}_{1}, \ldots, {\bf v}_{p})^{\top} $ is the corresponding orthogonal matrix with ${\bf v}_{i}$ being
the eigenvector associated to the $i$-th largest eigenvalue $\gamma_{i, p}, v_{i, 1} \geq 0, i=1, \ldots, p.$
Similarly, for the sample spectral decomposition, i.e., ${\bf S}={\bf U}{\bf L}{\bf U}^{\top}$, where ${\bf L}$ is a
diagonal matrix with eigenvalues $l_{i, p}, i=1, \ldots, p$, and ${\bf U}=({\bf u}_1, \ldots, {\bf u}_p)^{\top}$ is the
corresponding orthogonal matrix with ${\bf u}_{i}$ being the eigenvector corresponding  to
$l_{i, p}, u_{i1} \geq 0, i=1, \ldots, p$. Write ${\mbGamma}=\mbox{diag}(\gamma_{1, p}, \ldots, \gamma_{p, p})$,
and ${\bf L}=\mbox{diag}(l_{1, p}, \ldots, l_{p, p})$. Let $\Hat\mbPsi({\bf L})
=\mbox{diag}(\Hat\psi_{1}({\bf L}), \ldots, \Hat\psi_{p}({\bf L}))$ a real function of ${\bf L}$.

\vspace{0.3cm}
\indent {\bf Lemma}. (Perlman, 2007). {\it An estimator $\Hat{\mbSigma}$ is orthogonally equivariant if and only if}
\begin{align}
\Hat{\mbSigma}({\bf S})={\bf U}\Hat\mbPsi({\bf L}){\bf U}^{\top}.
\end{align}
\vspace{0.1cm}

\indent By the property of equality (2.4) and Perlman's comment that by restricting considerations to the class of
orthogonally equivariant estimators, the problem of estimating ${\mbSigma}$ reduces to that of estimating population
eigenvalues based on the sample eigenvalues (Perlman, 2007, page 67). Then, we may have the following:

\vspace{0.3cm}
\indent {\bf Proposition 2.1}. {\it Assume the loss function is convex, among the class of orthogonally equivariant
estimators then the minimum risk occurs at $R(\Hat{\mbSigma}, {\mbSigma})$ reduces to the minimum risk that of
$R(\Hat\mbPsi{_0}({\bf L}), {\mbGamma})$, where $\Hat\mbPsi{_0}({\bf L})$ is the estimator of ${\mbGamma}$ such that
$R(\Hat\mbPsi_{0}({\bf L}), {\mbGamma})$ is minimized}.
\vspace{0.2cm}

\indent Let $\Hat\mbPsi^{*}({\bf L})=\mbox{diag}(\Hat\psi^{*}_{1}({\bf L}), \ldots, \Hat\psi^{*}_{p}({\bf L}))$.
Among the class of orthogonally equivariant estimators, the risk function $R(\Hat {\mbSigma}, {\mbSigma})$ can then be
reduced to $R(\Hat{\mbPsi}^{*}({\bf L}),{\mbGamma})$. Consider the Stein loss function, the risk function is
\begin{align}
R(\Hat{\mbPsi}^{*}({\bf L}), {\mbGamma})& =\sum_{i=1}^{p}{\mathcal E}[\frac{\Hat{\psi}^{*}_{i}({\bf L})}
{\gamma_{i, p}}-\mbox{log}\frac{\Hat {\psi}^{*}_{i}({\bf L})}{{\gamma_{i, p}}}-1].
\end{align}
It is easy to see that the risk function (2.6) is minimized at ${\gamma}_{i, p}={\mathcal E}\{\Hat{\psi}^{*}_{i}({\b L})\},
i=1, \cdots, p$. When $p$ is fixed, it was shown that $l_{i, p}$ converges to $\gamma_{i, p}$ almost surely $(a.s.)$ as
$n\to \infty$ (Anderson, 2003). Hence, the best choice of $\Hat {\psi}^{*}_{i}({\bf L})$ among the class of orthogonally
equivariant estimators is to choose the consistent estimator $l_{i, p}$ of ${\gamma}_{i, p}, i=1, \ldots, p$. Namely, the
best choice of $\Hat \mbPsi{^*}({\bf L})$ is ${\bf L}$.

\vspace{0.3cm}
\indent {\bf Theorem 2.1}. {\it Among the class of orthogonally equivariant estimators, when the dimension $p$ is fixed.
Under the Stein loss function, the sample covariance matrix ${\bf S}$ is the best orthogonally equivariant estimator of
$\mbSigma$}.
\vspace{0.2cm}

\vspace{0.3cm}
\indent Stein (1975, 1986) proposed the well-known orthogonally equivariant nonlinear shrinkage estimators of the form
\begin{align}
\Hat {\mbSigma}_{S}&={\bf U}\Hat\mbPhi({\bf L}){\bf U}^{\top}, \mbox{where}~
\Hat\mbPhi({\bf L})=\mbox {diag}(\Hat\phi_{1}({\bf L}), \ldots, \Hat\phi_{p}({\bf L})) ~\mbox{with}  \\ \nonumber
 &\Hat\phi_{i}({\bf L})=nl_{i, p}(n-p+1-2l_{i, p}\sum_{j \ne i}\frac{1}{l_{j, p}-l_{i, p}})^{-1},
~ i=1, \ldots, p.
\end{align}
It is the golden standard, unfortunately, it has been pointed out that some of the $\Hat\phi_{i}({\bf L})$
might be negative and non-monotone numerically. To mitigate the problems, Stein recommended to use an isotonizing
algorithm procedure to adjust his estimators in (2.7).

\vspace{0.3cm}
\noindent {\bf 3. The Stein's estimator can be inadmissible when the dimension $p$ is fixed.}
\def \theequation{3.\arabic{equation}}
\setcounter{equation}{0}
\vspace{0.3cm}

\indent For the application of Proposition 2.1, we mainly adopt the Stein loss function (2.6).  Assumed that
$ \gamma_{1, p} > \ldots > \gamma_{p, p} > 0$, so that $l_{1, p} > \ldots > l_{p, p} > 0$ with probability one.

\indent First, we compare the best orthogonally equivariant estimator  of the population covariance matrix and the
Stein's estimator under Sten's loss function for the situation when the dimension $p$ is fixed. For simplicity,
we may assume that the population eigenvalues are widely dispersed, then the sample eigenvalues will also be
widely dispersed with probability one. Perlman (2007, page 70) pointed out the relationship that
$l_{i, p} \sum_{j \ne i}\frac{1}{l_{i, p}-l_{j, p}}\approx p-i$ holds when the population eigenvalues are widely
dispersed, $ i=1, \ldots, p$. In this situation, the Stein's estimator $\Hat {\mbSigma}_{S}$ in (2.7) is approximated
to the following form
\begin{align}
\Hat {\mbSigma}_{0}&={\bf U}\Hat\mbPhi^{0} ({\bf L}){\bf U}^{\top}, \mbox {where}~ \Hat\mbPhi^{0}({\bf L})
=\mbox {diag}(\Hat\phi^{0}_{1}({\bf L}), \ldots, \Hat\phi^{0}_{p}({\bf L})) ~\mbox {with}  \\ \nonumber
& \Hat\phi^{0}_{i}({\bf L})=nl_{i, p}(n+p-2i+1)^{-1}, ~ i=1, \ldots, p.
\end{align}

\indent After some algebraics the minimum risks of the best orthogonally equivariant estimator ${\bf S}$ and Stein's
estimator ${\Hat {\mbSigma}_{0}}$ become
\begin{align}
R_{m}({\bf L}, {\mbGamma})&=-\sum_{i=1}^{p}{\mathcal E}[\mbox{log}\frac{l_{i, p}}{{\gamma_{i, p}}}]  \\ \nonumber
 \mbox{and} \\
R_{m}(\Hat\mbPhi^{0}({\bf L}), {\mbGamma}) 
&=\sum_{i=1}^{p}[\frac{n}{n+p-2i+1}-\mbox{log}\frac{n}{n+p-2i+1}-1]+ R_{m}({\bf L}, {\mbGamma}),
\end{align}
respectively.

\indent Note that when $p$ is odd, then the middle term is located at $i=(p+1)/2$ and hence the middle term of
$\frac{n}{n+p-2i+1}-\mbox{log}\frac{n}{n+p-2i+1}-1, i=1, \ldots, p,$ becomes zero. After some straightforward
calculations, the minimum risk of $\Hat\mbPhi^{0}({\bf L})$ is
\begin{equation}
R_{m}(\Hat\mbPhi^{0}({\bf L}), {\mbGamma})=
   \sum_{i=1}^{[\frac{p}{2}]}[(2v_i-\mbox{log}v_i-2)]
  +R_{m}({\bf L}, {\mbGamma}), 
\end{equation}
where $[x]$ denotes the integer of $x$ and $v_i={1}/({1-(\frac{p-2i+1}{n})^{2}}), i=1, \ldots, [\frac{p}{2}]$. Note that
$v_i \geq 1, \forall i=1, \ldots, [\frac{p}{2}]$, and then $\sum_{i=1}^{[\frac{p}{2}]}(v_i-\mbox{log}v_i-1)\geq 0$. Thus,
\begin{align}
 R_{m}(\Hat\mbPhi^{0}({\bf L}), {\mbGamma})-R_{m}({\bf L}, {\mbGamma})
 & \geq 2\sum_{i=1}^{[\frac{p}{2}]}(v_i-\mbox{log}v_i-1)
   \\ \nonumber
 & \geq 0.
\end{align}
Therefore we have the following theorem.

\vspace{0.3cm}
\indent {\bf Theorem 3.1}. {\it Assume that the population eigenvalues are widely dispersed. Under the Stein loss
function when the dimension $p$ is fixed, then
$R_{m}(\Hat\mbPhi^{0}({\bf L}), {\mbGamma}) \geq R_{m}({\bf L}, {\mbGamma})$. And hence the Stein's estimator is
inadmissible.}
\vspace{0.2cm}

\vspace{0.3cm}
\noindent {\bf 4. High-dimensional case}
\def \theequation{4.\arabic{equation}}
\setcounter{equation}{0}
\vspace{0.3cm}

\indent For a large $(n, p)$ set-up, the large dimensional asymptotics framework is setup when $(n, p) \to \infty$
such that $c=\lim_{n \to \infty}p/n$ is fixed, $0 \leq c  < 1$. In this section, we extend the class of
orthogonally equivariant estimators to the realm of large dimensional asymptotics with the concentration $c \in (0,1)$
via the help of random matrix theory. Some basic notations of random matrix theory are presented in the following.

\vspace{0.3cm}
\noindent {\it 4.1 The Mar$\check{c}$enko-Pastur equation}
\vspace{0.2cm}

\indent The same as Ledoit and P$\acute{e}$ch$\acute{e}$ (2011), we make the following assumptions:

\indent A1. Note that ${\bf x}_{i}={\mbSigma}^{1/2}{\bf z}_{i}, i=1,  \ldots, n$, where ${\bf z}_{i}$ are independent
and identically distributed with mean $ {\bf 0}$ and covariance matrix ${\bf I}$. Assume that the 12th absolute central
moment of each variable $z_{ij}$ bounded by a constant.

\indent A2. The population covariance matrix $ {\mbSigma}$ is nonrandom positive definite.
$\mbox{lim inf}_{p \to \infty} \gamma_{p, p} > 0$ and $\mbox{lim sup}_{p \to \infty} \gamma_{1, p} < \infty$.

\indent A3. For large $(n,p)$ set-up, the large dimensional asymptotics framework is setted up when $(n,p) \to \infty$
such that $\lim_{n \to \infty} p/n=c \in (0, 1)$  in this paper.

\indent A4. Let $0 < \gamma_{p, p} < \cdots < \gamma_{1, p} $. The emperical spectral distribution of $\mbSigma$
defined by $H_{n}(\gamma)=\frac{1}{p}\sum_{i=1}^{p}1_{[\gamma_{i, p}, \infty)}(\gamma)$, converges as $p \to \infty$
to a probability distribution function $H(\gamma)$ at every point of continuity of $H$. The support of $H$,
$\mbox{Supp}(H)$, is included in a compact set $[h_1, h_2]$ with $0 < h_1 \leq h_2 < \infty$.

\indent Let $F_{n}(\lambda)=\frac{1}{p}\sum_{i=1}^{p}1_{[l_{i, p}, \infty)}(\lambda)$ be the sample spectral
distribution and $F$ be its limiting. Under the assumptions A1-A4, Mar$\check{c}$enko and Pastur (1967) proved that
$F_{n}$ converges to $F$ ~$a.s.$~ as $n \to \infty$.

\indent The Stieltjes transform of distribution function $F$ is defined by
\begin{align}
m_F(z)=\int_{-\infty}^{\infty}\frac{1}{l-z}dF(l), \forall z \in C^{+},
\end{align}
where $C^{+}$ is the half-plane of complex numbers with a strictly positive imaginary part. Let
\begin{align}
m_{F_{n}}(z)=p^{-1}\mbox{tr}[({\bf S}-z{\bf I})^{-1}],
\end{align}
then from the results of random matrix theory that $m_{F_{n}}(z)$ converges to $m_F(z)$ if and only if $F_n(z)$
converges to $F(z)$ weakly. Subsequently, the well-known Mar$\check{c}$enko-Pastur equation (Silverstein, 1995) in
literatures can be expressed in the following form
\begin{align}
m_F(z)=\int_{-\infty}^{\infty}\frac{1}{\gamma[1-c-czm_F(z)]-z}dH(\gamma), \forall z \in C^{+},
\end{align}
where $H$ denotes the limiting behavior of the population spectral distribution. Upon the Mar$\check{c}$enko-Pastur
equation, meaningful information of the population spectral distribution can be retrieved under the large
dimensional asymptotics framework. Choi and Silverstein (1995) further showed that
\begin{align}
\lim_{z\in C^{+}\to \l}m_F(z)=\check{m}_F(l)
\end{align}
exists for any $l \in R / \{0\}$.

\indent Let ${\Theta^{g}_{p}}(z)=p^{-1}\mbox{tr}[({\bf S}-z{\bf I})^{-1}g(\mbSigma)]$,
where $g(\cdot)$ is a scale function on the eigenvalues of a matrix such that
$g(\mbSigma)={\bf V}\mbox {diag}(g(\gamma_{1, p}), \ldots, g(\gamma_{p, p})){\bf V}^{\top}$. Ledoit and
P$\acute{e}$ch$\acute{e}$ (2011) proved that ${\Theta^{g}_{p}}(z)$ converges a.s. to ${\Theta^{g}}(z)$ under the
assumptions $A1-A4$ for any $z \in C^{+}$, where
\begin{align}
 {\Theta^{g}}(z)=\int_{-\infty}^{\infty}\frac{1}{\gamma[1-c-czm_F(z)]-z}g(\gamma)dH(\gamma),
    ~\forall z \in C^{+}.
\end{align}
Two special cases are $g(\mbSigma)={\mbSigma}$ and $g(\mbSigma)={\mbSigma}^{-1}$. Let
\begin{align}
\Delta_{p}(x)=\frac{1}{p}\sum_{i=1}^{p}{\bf u}^{\top}_{i}{\mbSigma}{\bf u}_{i}1_{\{l_i \leq x\}}, ~ x \in R/ \{0\}.
\end{align}
Ledoit and P$\acute{e}$ch$\acute{e}$ (2011) showed that $\Delta_{p}(x)$ converges a.s. to
$\Delta(x)=\int_{-\infty}^{x}\delta(l)dF(l)$  for $x \in  R / \{0\}$, where
$\delta(l)={l}/{|1-c-cl\check{m}(l)|^{2}}, ~l > 0, c \in (0, 1)$.
Also let
\begin{align}
\Delta^{(-1)}_{p}(x)=\frac{1}{p}\sum_{i=1}^{p}{\bf u}^{\top}_{i}{\mbSigma}^{-1}{\bf u}_{i}1_{\{l_i \leq x\}},
~ x \in R/ \{0\}.
\end{align}
Ledoit and P$\acute{e}$ch$\acute{e}$ (2011) also showed that $\Delta^{(-1)}_{p}(x)$ converges a.s. to
$\Delta^{(-1)}(x)=\int_{-\infty}^{x}\delta^{(-1)}(l)dF(l)$ for $x \in  R / \{0\}$, where
$\delta^{(-1)}(l)=({1-c-2cl\mbox{Re}[\check{m}(l)]})/{l}, ~l > 0, c \in (0, 1)$. They also claimed that
 $d^{*}_i={\bf u}^{\top}_{i}{\mbSigma}{\bf u}_{i}$ and $a^{*}_{i}={\bf u}^{\top}_{i}{\mbSigma}^{-1}{\bf u}_{i}$
can be approximated by ${l_{i}}/{|1-c-cl_i\check{m}(l_i)|^{2}}$ and
$({1-c-2cl_{i}\mbox{Re}[\check{m}(l_{i})]})/{l_{i}}, i=1, \ldots, p$, respectively, where $l_{i}$
denotes the $(1-\alpha)$-quantile of limiting sample spectral distribution $F$ with $\alpha \in (0, 1)$,
so that $[p(1-\alpha)]=i$, with $[x]$ denoting the largest integer of $x$. Note that
$(\delta(l))^{-1} \ne \delta^{(-1)}(l)$, which is against the statistical sense. For the consistency of notations,
we may write that
\begin{align}
\phi_{i}({\bf L})=\frac{1}{\delta^{(-1)}(l_{i})}=\frac{l_{i}}{1-c-2cl_{i}\mbox{Re}[\check{m}(l_{i})]}, ~i=1, \ldots, p.
\end{align}
Ledoit and Wolf (2018) neither treated $\delta(l_{i})$ as the parameter $\gamma_{i}$ defined in (4.18), nor treated
$\delta^{(-1)}(l_{i})$ as the parameter $\gamma^{-1}_{i}, i=1, \ldots, p$. Instead, they and Stein (1975, 1986) all
treated $\phi_{i}({\bf L}), i=1, \ldots, p$, as the interesting parameters to be estimated.

\indent Stein's estimator $\Hat\phi_{i}({\bf L})$ in (2.7) is the consistent estimator of $\phi_{i}({\bf L})$ in (4.8),
$i=1, \ldots, p.$ Ledoit and Wolf (2018) commented that Stein's estimator has theoretical limitations, and hence they
recommended to adopt their proposed estimator of ${\mbSigma}$, which is of the form
\begin{align}
\Hat{\mbSigma}_{LW}&={\bf U}\Hat\mbPhi^{*}({\bf L}){\bf U}^{\top},\mbox {where}~\Hat\mbPhi^{*}({\bf L})
=\mbox {diag}(\Hat\phi^{*}_{1}({\bf L}), \ldots, \Hat\phi^{*}_{p}({\bf L})) ~\mbox{with} \\ \nonumber
 &\Hat\phi^{*}_{i}({\bf L})=\frac{l_{i, p}}{1-\frac{p}{n}-2\frac{p}{n}l_{i, p}
\mbox{Re}[\check{m}^{\Hat \mbtau_{n}}_{n, p}(l_{i, p})]}, ~i=1, \ldots, p,
\end{align}
where $\mbox{Re}[\check{m}^{\Hat \mbtau_{n}}_{n, p}(l_{i, p})]$ is the consistent estimator of
$\mbox{Re}[\check{m}_F(l_{i})]$, as well as hinges on a multivariate quantized eigenvalues sample function. Readers
may see their paper for details (Ledoit and Wolf, 2012). They claimed their estimator in (4.9) naturally circumvents
the disadvantages of Stein's estimator in (2.7), and their estimator performs better than that of Stein's estimator
in practice via the evidence of extensive Monte-Carlo simulations. They pointed out that both the estimators in (2.7)
and in (4.9) are different forms to estimate the Cauchy principal value $\mbox{Re}[\check{m}(l_{i})]$ in (4.8).
The estimator of Ledoit and Wolf in (4.9) uses the product of cutting-edge research of random matrix theory,
while Stein uses the naive empirical distribution of sample eigenvalues. Since the Mar$\check{c}$enko-Pastur distribution
function is continuous, it is quite reasonable to expect that the performance of using a smoothed version estimator to
estimate the Hilbert transform of the distribution function will be better in practice. Theoretically, both two different
estimators should enjoy the same large-dimensional asymptotics optimality because both $\Hat\phi^{*}_{i}({\bf L})$ and
$\Hat\phi_{i}({\bf L})$ are all claimed to be the consistent estimators of $\phi_{i}({\bf L})$ defined in (4.8),
$i=1, \ldots, p$. However, we may concern that $\phi_{i}({\bf L})$ might not be exactly the same as the eigenvalues
$\gamma_{i, p}$ of the population covariance matrix, $\forall ~i=1, \ldots, p$.

\indent When both the dimension $p$ and the sample size $n$ are large so that $\lim_{n \to \infty}p/n=c \in (0, 1)$,
then ${\bf L}$ is no longer to be the consistent estimator of ${\mbGamma}$. The oracle parameters
$d^{*}_{i}={\bf u}^{\top}_{i}{\mbSigma}{\bf u}_{i}$ and
$a^{*}_{i}={\bf u}^{\top}_{i}{\mbSigma}^{-1}{\bf u}_{i}, i=1, \ldots, p$, which proposed to be estimated via the
projection method by Ledoit and P$\acute{e}$ch$\acute{e}$ (2011), are generally not the functions of the population
eigenvalues only. Hence we may notice two things. First, both $d^{*}_i$ and $a^{*}_i$ are the functions that are also
involved with both the sample and the population eigenvectors, hence Proposition 2.1 is not applicable for both of
them unless ${\bf U} \to {\bf V}~ a.s.$, which remains unsolved when the dimension $p$ is large in the literature. Secondly,
even if it is to be true, it is quite reasonable to expect that $\delta(l_{i})\delta^{(-1)} (l_{i})=1, i=1, \ldots, p,$ under
the large dimensional asymptotics setup, however, this anticipation could not be accomplished via the results of Ledoit
and Wolf (2012). Thus, people might suspect that the approach of Ledoit and Wolf can not be consistent with the result of
Proposition 2.1. As such, it is still an unsolved problem whether the consistent estimators of eigenvalues ${\gamma}_{i, p}$
of the population covariance matrix exist, $i=1, \ldots, p.$ In the next section, we further study whether another kind of
type $\phi_{i}({\bf L})$ exists so that it is equivalent to the eigenvalues $\gamma_{i}$, $i=1, \ldots, p$ under the large
dimensional asymptotics setup. We provide an affirmative answer by looking into the Mar$\check{c}$enko-Pastur equation
in the next subsection.

\vspace{0.3cm}
\noindent {\it 4.2 The identity equation between quantiles $l$ and $\gamma$.}
\vspace{0.2cm}

\indent We project the sample eigenvectors onto the population eigenvectors the same as Ledoit and P$\acute{e}$ch$\acute{e}$
(2011) did and adopt their notations, and then study the functionals of the type, $\forall z \in C^{+}$,
\begin{align}
{\Theta}^{1_{(-\infty, \gamma)}}_{p}(z)=\frac{1}{p}\sum_{i=1}^{p}\frac{1}{l_{i, p}-z}\sum_{j=1}^{p}
 |{\bf u}^{\top}_{i} {\bf v}_{j}|^{2}1_{(-\infty, \gamma)}.
\end{align}
Thus, ${\Theta}^{1_{(-\infty, \gamma)}}_{p}(z)$ converges to ${\Theta}^{1_{(-\infty, \gamma)}}(z)$ a.s. as $p \to
\infty$,
$\forall z \in C^{+}$, where
\begin{align}
{\Theta}^{1_{(-\infty, \gamma)}}(z)=\int_{-\infty}^{\gamma}\frac{1}{t[1-c-czm_F(z)]-z}dH(t).
\end{align}

\indent Consider a bivariate distribution function
\begin{align}
\Phi_{p}(\lambda, \gamma)=\frac{1}{p}\sum_{i=1}^{p}\sum_{j=1}^{p}
 |{\bf u}^{\top}_{i} {\bf v}_{j}|^{2}1_{[\lambda_{i}, \infty)}(\lambda)1_{[\gamma_{i}, \infty)}(\gamma),
\end{align}
$(\lambda, \gamma) \in R^{2}$. Thus,
\begin{align}
\Phi_{p}(\lambda, \gamma)=\lim_{\eta \to 0^{+}}\int_{-\infty}^{\lambda}\frac{1}{\pi}
\mbox{Im}[{\Theta}^{1_{(-\infty, \lambda)}}_{p}(l+i\eta)]dl
\end{align}
holds. Therefore, $\lim_{p \to \infty} \Phi_{p}(\lambda, \gamma)$ exists and is equal to
\begin{align}
\Phi(\lambda, \gamma)=\lim_{\eta \to 0^{+}}\int_{-\infty}^{\lambda}\frac{1}{\pi}
\mbox{Im}[{\Theta}^{1_{(-\infty, \lambda)}}(l+i\eta)]dl,
\end{align}
for evry $(\lambda, \gamma)\in R^{2}$ where $\Phi$ is continuous.

\indent Let $a=\mbox{Re}[1-c-cl\check{m}_F(l)]$ and $b=\mbox{Im}[1-c-cl\check{m}_F(l)]$. Since
$F'(l)=\frac{1}{\pi}\mbox{Im}[\check{m}_F(l)]$, thus $b=-\pi clF'(l)$. When $c \in (0, 1)$, Ledoit and
P$\acute{e}$ch$\acute{e}$ (2011) proved that
\begin{align}
\Phi(\lambda, \gamma)=\int_{-\infty}^{\gamma}\int_{-\infty}^{\lambda}\frac{clt}{(at-l)^{2}+b^{2}t^{2}}dF(l)dH(t),
\end{align}
which is the bivariate distribution function of $\lambda$ and $\gamma$.

\indent In order that integrals evaluated over the entire real line are convergent, we may assume that $F'(l) \to 0 $
sufficiently fast as $|l| \to \infty ~$. Thus, the marginal distribution of $\gamma$, by putting $\lambda \to \infty$
in $\Phi(\lambda, \gamma)$, is
\begin{align}
\Phi(\gamma)=\Phi(\infty, \gamma)&=\int_{-\infty}^{\gamma}\int_{-\infty}^{\infty}
\frac{clt}{(at-l)^{2}+b^{2}t^{2}}dF(l)dH(t)\\ \nonumber
&=\int_{-\infty}^{\gamma}\int_{-\infty}^{\infty}\frac{1}{\pi}
\frac{-bt}{(at-l)^{2}+b^{2}t^{2}}dldH(t).
\end{align}
Under the regularity conditions that $F^{'}(l)$ vanishes as $|l| \to \infty$. As $\lambda \to \infty$, then the
infinitestimal quantity $F'(l) \to 0$, which implies that $-bt \to 0$ for any finite $t$. Write $\epsilon=-bt$.
Sending $\epsilon \to 0$ and mimicking the proof of the imaginary part of Sokhotski-Plemelj formula, with the
representation in terms of the delta function
${\delta}(x)=\frac{1}{\pi}\lim_{\epsilon \to 0}\frac{\epsilon}{x^{2}+{\epsilon}^{2}}$. Then, we have
\begin{align}
\Phi(\gamma)&= \int_{-\infty}^{\gamma}\int_{-\infty}^{\infty}\frac{1}{\pi}
\lim_{\epsilon \to 0}\frac{\epsilon}{(at-l)^{2}+{\epsilon}^{2}}dldH(t) \\ \nonumber
&=\int_{-\infty}^{\gamma}\int_{-\infty}^{\infty}\delta(at-l)dldH(t) \\ \nonumber
&=H(\gamma), \mbox{with the equality equation that}~ \gamma=\frac{l}{a} ~\mbox{being held}.
\end{align}
We prove the marginal distribution of $\gamma$ is $H(\gamma)$, and accompany with the equality equation. Hence,
we may establish the equality equation of the quantiles in the following.

\vspace{0.3cm}
\indent {\bf Theorem 4.1}. {\it  Let $\gamma_{i}$ and $l_{i}$ denote the $(1-\alpha)$-quantiles of limiting population
and sample spectral distributions $H$ and $F$, respectively, $\alpha \in (0, 1)$, so that $[p(1-\alpha)]=i$, with $[x]$
denoting the largest integer of $x$, and $\mbox{Re}[\check{m}_F(l_{i})]$ denotes the Cauchy principle value of
$\check{m}_{F}(l_{i})$. Under the assumptions A1-A4. Assume that $F^{'}(l)$ vanishes as $|l| \to \infty$. When
$\lim_{n \to \infty} p/n=c \in (0, 1)$, then we have
\begin{align}
\gamma_{i}= \frac{l_{i}}{1-c-cl_{i}[\mbox{Re}\check{m}_F(l_{i})]}, i=1, \ldots, p.
\end{align}
}

\indent Note that $\gamma_{i}$ in (4.18) is not exactly the same as ${\phi}_{i}(\bf L)$ defined in (4.8), $i=1, \cdots, p.$

\vspace{0.3cm}
\noindent {\it 4.3 The consistent estimators of population eigenvalues}
\vspace{0.2cm}

\indent  Next, the main work is to estimate the population eigenvalues $\gamma_{i, p}$ based on the sample eigenvalues
$l_{i, p},  i=1,\ldots, p$. Let $\mbox{Supp}(F)$ be the support of $F$. Via Theorem 4.1 of Choi and Silverstein (1995),
Ledoit and P$\acute{e}$ch$\acute{e}$ (2011) pointed out that if $l_{i} \notin \mbox{Supp}(F)$, then
${l_{i}}/{1-c-cl_{i}[\mbox{Re}\check{m}_F(l_{i})]} \notin \mbox{Supp}(H),$ for $l_{i} \in R/\{0\}, i=1, \ldots, p$.

\indent  Write $\gamma_{i}=\psi_{i}({\bf L}), i=1, \ldots, p$, based on the results of (4.18), we may then propose
another kind of orthogonally equivariant estimator $\Hat {\mbSigma}_{T}$ of $\mbSigma$, which is of the form
\begin{align}
\Hat {\mbSigma}_{T}={\bf U}\Hat\mbPsi({\bf L}){\bf U}^{\top}, &~ \mbox {where}~
\Hat\mbPsi({\bf L})=\mbox {diag}(\Hat\psi_{1}({\bf L}), \cdots, \Hat\psi_{p}({\bf L})) ~\mbox {with}  \\ \nonumber
\Hat\psi_{i}({\bf L})&=\frac{nl_{i, p}}{n-p+1-pl_{i, p}{\check{m}}_{F_n}(l_{i, p})}  \\ \nonumber
         &=nl_{i, p}(n-p+1-l_{i, p}\sum_{j \ne i}\frac{1}{l_{j, p}-l_{i, p}})^{-1}, ~ i=1, \ldots, p.
\end{align}

\indent Note that the estimator in (4.19) is very similar to Stein's estimator in (2.7). The consistent property of
$\Hat\psi_{i}({\bf L})$ can then be followed by the same arguments as those of the Stein's estimator. Under the assumptions
A1- A4, Mar$\check{c}$enko and Pastur (1967) proved that $F_{n}(x)$ converges to $F(x)$ a.s., thus $l_{i, p}$ converges
to $l_{i}~ a.s., i=1, \ldots, p$ and $m_{F_{n}}(z)$ converges to $m_F(z)~ a.s.$ when $p$ is large. And hence by (4.2), we
have that ${\check{m}}_{F_n}(l_{i, p})=\frac{1}{p}\sum_{j \ne i}\frac{1}{l_{j, p}-l_{i, p}}$. Note that $l_{i, p}$ is the
consistent estimator of $l_{i}$, $ i=1, \ldots, p$. And then ${\check{m}}_{F_n}(l_{i, p})$ is the consistent estimator
of $\mbox{Re}[\check{m}_F(l_{i})], i=1, \ldots, p$. Thus, $\Hat\psi_{i}({\bf L})$ is the consistent estimator of
$\psi_{i}({\bf L})~(i.e., \gamma_{i}), i=1, \ldots, p.$ By the assumption A4 that $H_n$ converges to $H$, thus
$\gamma_{i, p}$ converges to $\gamma_{i},  i=1, \ldots, p$, when $p$ is large. To estimate $\gamma_{i,p}$ can be
viewed the same as to estimate $\gamma_{i}$ when $p$ is large, $i=1, \ldots, p$. Namely, $\Hat\psi_{i}({\bf L})$ can be
viewed the same as the consistent estimator of ${\gamma}_{i, p}, i=1, \ldots, p.$ Therefore, we have the following.

 \vspace{0.3cm}
\indent {\bf Theorem 4.2}. {\it Under the assumptions of Theorem 4.1.  Let
${\mbGamma}^{0}=\mbox {diag}(\gamma_{1}, \ldots, \gamma_{p})$ and
$\Hat\mbPsi({\bf L})=\mbox {diag}(\Hat\psi_{1}({\bf L}), \ldots, \Hat\psi_{p}({\bf L}))$ be defined in (4.18) and
(4.19), respectively. When $\lim_{n \to \infty} p/n=c \in (0, 1)$, then $\Hat\mbPsi({\bf L})$ is the consistent estimator
of ${\mbGamma}^{0}$, namely $\Hat\mbPsi({\bf L})$ is the consistent estimator of ${\mbGamma}$}.
\vspace{0.3cm}

\indent {\bf Remarks}. The results of Theorem 4.2 make up the deficiency of Stein's estimator $\Hat\mbPhi({\bf L})$ in
(2.7) and Ledoit and Wolf's estimator $\Hat\mbPhi^{*}({\bf L})$ in (4.9), which are not consistent for $\mbPsi({\bf L})$.
\vspace{0.3cm}

\indent Consider the normalized Stein loss function (i.e., the one in (2.6) divided by $p$). Then the risk function of
$\Hat {\mbSigma}_{T}$ is
\begin{align}
R^{*}(\Hat\mbPsi({\bf L}), {\mbGamma})=\frac{1}{p}\sum_{i=1}^{p}{\mathcal E}[\frac{nl_{i, p}}
{(n-p+1-pl_{i, p}{\check{m}}_{F_n}(l_{i, p})){\gamma}_{i, p}}-\mbox{log}\frac{nl_{i, p}}
{(n-p+1-pl_{i, p}{\check{m}}_{F_n}(l_{i, p})){\gamma}_{i, p}}-1].
\end {align}

\indent When $p$ is large such that $p/n \to c\in (0, 1)$, by Proposition 2.1 we may note that the risk function is
minimized at $\gamma_{i, p}={\mathcal E}\{\frac{nl_{i, p}}{n-p+1-pl_{i, p}\check{m}_{F_{n}}(l_{i, p})}\}, i=1, \ldots, p$.
Since $\Hat\psi_{i}({\bf L})$ is the consistent estimator of ${\gamma}_{i}$ and  ${\gamma}_{i, p} \to {\gamma}_{i}$ by
the assumption A4, hence the best choice of $\Hat \psi^{*}_{i}({\bf L})$ among the class of orthogonally equivariant
estimators is $\Hat \psi_{i}({\bf L}), i=1, \ldots, p$. Thus, we have the following.

\vspace{0.3cm}
\indent {\bf Theorem 4.3}. {\it Under the assumptions of Theorem 4.1. When $\lim_{n \to \infty} p/n=c \in (0, 1)$,
under the Stein loss function with being normalized by $1/p$,  the estimator $\Hat {\mbSigma}_{T}$ is the best
orthogonally equivariant estimator of the population covariance matrix}.
\vspace{0.2cm}

\vspace{0.3cm}
\noindent {\it 4.4 Dominations}
\vspace{0.2cm}

\indent Note that the Ledoit and Wolf's estimator $\Hat {\mbSigma}_{LW}$ in (4.9) is asymptotically equivalent to the
Stein's estimator in (2.7) when both $p$ and $n$ are large so that $\lim_{n \to \infty}p/n=c \in (0, 1)$. Under the
Stein's loss function with the normalization by $1/p$ when $p/n \to c \in (0, 1)$, with $l_{i, p}\sum_{j \ne i}
\frac{1}{l_{i, p}-l_{j, p}}$ being replaced by $p-i, ~\forall i=1, \ldots, p$, then after some straightforward
algebraic the minimum risks of Stein estimator $\Hat {\mbSigma}_{0}$ in (3.1) and sample covariance matrix ${\bf S}$
are of the forms
\begin{align}
R^{*}_{m}(\Hat\mbPhi^{0}({\bf L}), {\mbGamma})= &
\frac{1}{p}[\sum_{i=1}^{p}\frac{n-i+1}{n+p-2i+1}+\sum_{i=1}^{p}
   \mbox{log}\frac {n+p-2i+1}{n}-p-\sum_{i=1}^{p}{\mathcal E}[\mbox{log}\frac{l_{i, p}}{{\gamma_{i, p}}}]] \\ \nonumber
&\mbox{and} \\
R^{*}_{m}({\bf L}, {\mbGamma})
&=\frac{1}{p}[\sum_{i=1}^{p}\frac{n-i+1}{n}-p
   -\sum_{i=1}^{p}{\mathcal E}[\mbox{log}\frac{l_{i, p}}{{\gamma_{i, p}}}]] \ \nonumber
\end{align}
respectively.

\indent Similarly, after some straightforward calculations the minimum risk of $\Hat {\mbSigma}_{T}$ is
\begin{align}
R^{*}_{m}(\Hat\mbPsi({\bf L}), {\mbGamma})=\frac{1}{p}[\sum_{i=1}^{p}\mbox{log}\frac{n-i+1}{n}-\sum_{i=1}^{p}
{\mathcal E}[\mbox{log}\frac{l_{i, p}}{{\gamma_{i, p}}}]],
\end{align}

\indent Let $x_i=({n-i+1})/{n+p-2i+1}$ and $y_i=({n-i+1})/n, i=1, \ldots, p$. Note that
$0 < x_i \leq 1$ and $0 < y_i \leq 1, i=1, \ldots, p$, and hence $x_i-\mbox{log}x_i-1 \geq 0$ and
$y_i-\mbox{log}y_i-1 \geq 0, i=1, \ldots, p$. Thus we have
\begin{align} \nonumber
R^{*}_{m}(\Hat\mbPhi^{0}({\bf L}), {\mbGamma})& =\frac{1}{p}[\sum_{i=1}^{p}(x_i-\mbox{log}x_i-1)]
    +R^{*}_{m}(\Hat\mbPsi({\bf L}), {\mbGamma}) \\ \nonumber
    & \geq R^{*}_{m}(\Hat\mbPsi({\bf L}), {\mbGamma}),  \\ \nonumber
\hspace{0.3cm} \mbox{and}   \\ \nonumber
R^{*}_{m}({\bf L}, {\mbGamma})&= \frac{1}{p}[\sum_{i=1}^{p}(y_i-\mbox{log}y_i-1)]
    +R^{*}_{m}(\Hat\mbPsi({\bf L}), {\mbGamma}) \\ \nonumber
  & \geq R^{*}_{m}(\Hat\mbPsi({\bf L}), {\mbGamma}),
\end{align}
respectively. Hence, we have the following theorem.

\vspace{0.3cm}
\indent {\bf Theorem 4.4}. {\it Under the assumptions of Theorem 4.1. Assume that the population eigenvalues are widely
dispersed. Under the normalized Stein loss function when both $p$ and $n$ are large so that
$\lim_{n \to \infty}p/n=c \in (0, 1)$, then
$R^{*}_{m}(\Hat\mbPsi({\bf L}), {\mbGamma}) \leq R^{*}_{m}(\Hat\mbPhi^{0}({\bf L}), {\mbGamma})$ and
$R^{*}_{m}(\Hat\mbPsi({\bf L}), {\mbGamma}) \leq R^{*}_{m}({\bf L}, {\mbGamma})$, respectively. Both Stein's estimator
$\Hat {\mbSigma}_{0}$ and the sample covariance matrix ${\bf S}$ are inadmissible.}
\vspace{0.2cm}

\vspace{0.3cm}
\noindent {\bf Acknowledgements.}
\vspace{0.2cm}

\vspace{0.3cm}
\noindent {\bf References}
\vspace{0.2cm}

\begin{enumerate}
\item Anderson, T. W. (2003).  {\it An Introduction to Multivariate Statistical Analysis,} 3rd edition,
Wiely, New York.


\item Choi, S. I. and Silverstein, J. W. (1995). Analysis of the limiting spectral distribution of large
dimensional random matrices. J. Multivariate Anal. 54, 295-309.

\item James, W. and Stein, C. (1961). Estimation with quadratic loss. Proc. Fourth Berkeley Symp. Math.
Statist. Probab. 1, 361-379. California Press, Berkeley, CA.

\item Ledoit, O. and P$\acute{e}$ch$\acute{e}$, S. (2011). Eigenvectors of some large sample covariance matrix
ensembles. Probab. Theory Relat. Fields 151, 233-264.

\item Ledoit, O. and Wolf, M. (2012). Nonlinear shrinkage estimation of  large-dimensional covariance matrices.
Ann. Statist. 40, 1024-1060.

\item Ledoit, O. and Wolf, M. (2018). Optimal estimation of a large-dimensional covariance matrix under Stein's loss.
  Bernoulli 24, 3791-3832.

\item Mar$\check{c}$enko, V. A. and Pastur, L. A. (1967). Distribution of eigenvalues for some sets of random
matrices. Sb. Math. 1, 457-483.

\item Perlman, M. D. (2007). {\it Multivariate Statistical Analysis,} Univ. Washington, Seattle, Washington.

\item Pourahmadi, Mohsen. (2013). {\it High-Dimensional Covariance Estimation}. Wiley, New York.


\item Rajaratnam, B. and Vincenzi, D. (2016). A theoretical study of Stein's covariance estimator.
 Biometrika 103, 653-666 (2016).

\item Silverstein, J. W. (1995). Strong convergence of the emperical distribution of eigenvalues of large dimensional
random matrices. J. Multivariate Anal. 55, 331-339.

\item Stein, C. (1956). Inadmissibility of the usual estimator of the mean of a multivariate normal distribution.
Proc. Third Berkeley Symp. Math. Statist. Probab. 1, 197-206. California Press, Berkeley, CA.

\item Stein, C. (1975). Estimation of a covariance matrix. Rietz lecture, 39th Annual Meeting IMS,
 Atalanta, Georgia.

\item Stein, C. (1986). Lectures on the theory of estimation of many parameters. J. Math. Sci. 43, 1373-1403.


\item Zagidullina, A. (2021). {\it High-Dimensional Covariance Matrix Estimation: An Introduction to Random Matrix Theory}.
  SpringerBriefs in Applied Statistics and Econometrics. Switzerland.

\end{enumerate}

\end{document}